\documentclass[11pt]{article}
\usepackage{amsfonts,amsmath,amssymb,amssymb,tikz}
\usetikzlibrary{arrows}
\newtheorem{theorem}{\bf Theorem}[section]
\newtheorem{corollary}[theorem]{\bf Corollary}
\newtheorem{lemma}[theorem]{\bf Lemma}
\newtheorem{proposition}[theorem]{\bf Proposition}

\newcommand{\proof}{\noindent{\bf Proof.\ }}
\newcommand{\qed}{\hfill $\square$ \bigskip}
\newcommand{\itemqed}{\hfill $\square$}
\newcommand{\mathqed}{\qquad \square}

\newcommand{\divides}{\,|\,}
\newcommand{\notdiv}{\!\not\,\!|\ }

\begin{document}

\title{\bf Vertex and edge orbits of Fibonacci and Lucas cubes}

\author{
Ali Reza Ashrafi $^{a}$
\and
Jernej Azarija $^{b}$
\and
Khadijeh Fathalikhani $^{a}$
\and
Sandi Klav\v zar $^{c,d,b}$
\and 
Marko Petkov\v sek $^{c,b}$
}

\date{}

\maketitle

\begin{center}
$^a$ Department of Pure Mathematics, Faculty of Mathematical  \\
Sciences, University of Kashan, Kashan, Iran \\
\medskip

$^b$ Institute of Mathematics, Physics and Mechanics, Ljubljana, Slovenia\\
\medskip

$^c$ Faculty of Mathematics and Physics, University of Ljubljana, Slovenia\\
\medskip

$^d$ Faculty of Natural Sciences and Mathematics, University of Maribor, Slovenia
\end{center}

\begin{abstract}
The Fibonacci cube $\Gamma_n$ is obtained from the $n$-cube $Q_n$ by removing all the vertices that contain two consecutive 1s. If, in addition, the vertices that start and end with 1 are removed, the Lucas cube $\Lambda_n$ is obtained. The number of vertex and edge orbits, the sets of the sizes of the orbits, and the number of orbits of each size, are determined for the Fibonacci cubes and the Lucas cubes under the action of the automorphism group. In particular, the set of the sizes of the vertex orbits of $\Lambda_n$ is $\{k \ge 1;\ k \divides n\} \cup\, \{k \ge 18;\ k \divides 2n\}$, 
the number of the vertex orbits of $\Lambda_n$ of size $k$, where $k$ is odd and divides $n$, is equal to $\sum_{d\divides k}\mu\left(\frac{k}{d}\right) F_{\lfloor \frac{d}{2}\rfloor + 2}$, 
and the number of the edge orbits of $\Lambda_n$ is equal to the number of the vertex orbits of $\Gamma_{n-3}$ when $n \ge 5$. Primitive strings, dihedral transformations and asymmetric strings are essential tools to prove these results.        
\end{abstract}


\medskip\noindent
{\bf AMS Subj. Class:} 68R15, 05C30

%%%%%%%%%%%%%%%%%%%%%%%%%%%%%%%%%%%%%%%%%%%%%%%%%%%%%%%%%%%%%%%%%%%%%
%%%%%%%%%%%%%%%%%%%%%%%%%%%%%%%%%%%%%%%%%%%%%%%%%%%%%%%%%%%%%%%%%%%%%
\section{Introduction}
%%%%%%%%%%%%%%%%%%%%%%%%%%%%%%%%%%%%%%%%%%%%%%%%%%%%%%%%%%%%%%%%%%%%%
%%%%%%%%%%%%%%%%%%%%%%%%%%%%%%%%%%%%%%%%%%%%%%%%%%%%%%%%%%%%%%%%%%%%%

Fibonacci cubes $\Gamma_n$~\cite{hsu-93} and the closely related Lucas cubes $\Lambda_n$~\cite{mupe-2001} have been investigated from many points of view, let us briefly overview some recent achievements. Formulas for the number of vertices of a given degree as well as the corresponding generating functions were determined in~\cite{klmope-2011}, while the domination number and the 2-packing number of these cubes were studied in~\cite{cakl-2011,pizo-2012}. Motivated by the structure of $\Gamma_n$ and $\Lambda_n$ as interconnection networks, Mollard~\cite{mo-2012} characterized maximal induced hypercubes in these cubes and also determined the number of such hypercubes. From the metric graph theory point of view, eccentricity sequences were obtained in~\cite{camo-2012}, the Wiener index and the Hosoya polynomial were determined in~\cite{klmo-2012}, while in~\cite{klmo-2014} the asymptotic average eccentricity was determined. In the latter paper it is also proved that the eccentricity of a vertex of a given Fibonacci cube is equal to the depth of the associated leaf in the corresponding Fibonacci tree. For a connection between Fibonacci/Lucas cubes and and Hasse diagrams (of the independent subsets of powers of paths and cycles) see~\cite{coot-2013}. From the chemical graph theory perspective we point out that Lucas cubes turned out to be precisely the so-called  resonance graphs of cyclic fibonacenes~\cite{zibe-2013a}. Very recently several advances were made also from the algorithmic point of view. Linear recognition algorithm for Fibonacci cubes and for Lucas cubes were developed by Vesel~\cite{ve-2014} and Taranenko~\cite{ta-2013}, respectively, while Ramras~\cite{ra-2014} studied off-line routing of linear permutations on these cubes. For additional  information on Fibonacci cubes see the survey~\cite{kl-2013}. 

There are several reasons for this wide interest. These cubes are induced sugraphs of hypercubes that inherit many of the fine properties of the latter class. The main tool to derive such properties for Fibonacci cubes is the so-called fundamental decomposition that decomposes $\Gamma_n$
into $\Gamma_{n-1}$ and $\Gamma_{n-2}$, similarly as the $n$-cube decomposes into two $(n-1)$-cubes via the Cartesian product operation. 
(There is also a similar decomposition for Lucas cubes.) On the other hand, the order of Fibonacci/Lucas cubes grows much slower than that
of hypercubes, a property important for interconnection networks. A strong source of interest for these cubes also comes from theoretical chemistry,
where Fibonacci cubes are precisely the so-called resonance graphs of fibonaccenes~\cite{klzi-2005} (see~\cite{zhou-2009} for a generalization of this
result), while for the role of Lucas cubes in chemistry, besides the already mentioned paper~\cite{zibe-2013b}, see also~\cite{zibe-2012,zibe-2013b}. Fibonacci cubes also led to the notion of the Fibonacci dimension of a graph~\cite{caep-2011,ve-2013}. 

When it comes to symmetries, it seems that only the automorphism groups of Fibonacci and Lucas cubes have been determined so far~\cite{cakl-2011}. Hence, in this paper we take a closer look at their symmetries, more precisely at the orbits under the action of the automorphism group. We proceed as follows. The next two sections are of preliminary nature. In the first of them we introduce concepts and notations needed, and recall or prove some related results. In the subsequent section we investigate some properties of dihedral transformations of nonempty strings defined over a finite alphabet. In Section~\ref{sec:orbit-fibonacci}, we determine the number of vertex and edge orbits, the sets of the sizes of the orbits, and the number of orbits of each size of Fibonacci cubes, as well as give a combinatorial interpretation for the number of vertex orbits. In the last section we prove parallel results for Lucas cubes. Contrary to Fibonacci cubes, where there are only orbits of size 1 and 2, the situation with Lucas cubes is more intriguing and complex.

%%%%%%%%%%%%%%%%%%%%%%%%%%%%%%%%%%%%%%%%%%%%%%%%%%%%%%%%%%%%%%%%%%%%%
%%%%%%%%%%%%%%%%%%%%%%%%%%%%%%%%%%%%%%%%%%%%%%%%%%%%%%%%%%%%%%%%%%%%%
\section{Preliminaries}
%%%%%%%%%%%%%%%%%%%%%%%%%%%%%%%%%%%%%%%%%%%%%%%%%%%%%%%%%%%%%%%%%%%%%
%%%%%%%%%%%%%%%%%%%%%%%%%%%%%%%%%%%%%%%%%%%%%%%%%%%%%%%%%%%%%%%%%%%%%

In this section we formally introduce the cubes studied here, list some notation, and prove an identity on the Lucas numbers. The $n$-cube $Q_n$, $n\ge 0$, is the graph whose
vertex set is the set of all binary strings of length $n$, two vertices being adjacent if they differ in exactly one position. The $n$-dimensional
{\em Fibonacci cube\/} $\Gamma_n$ is the subgraph of $Q_n$ induced by the set of all vertices that have no two consecutive 1s. Strings with no two
consecutive 1s are called {\em Fibonacci strings}. The $n$-dimensional {\em Lucas cube\/} $\Lambda_n$ is obtained from $\Gamma_n$ by removing all the vertices 
that begin and end with 1. The vertices of Lucas cubes are called {\em Lucas strings}. The Fibonacci numbers $F_n$ are defined by $F_0 = 0$, $F_1 = 1$,
and $F_n = F_{n-1} + F_{n-2}$, $n\ge 2$, and the Lucas numbers $L_n$ by $L_0 = 2$, $L_1 = 1$, and $L_n = L_{n-1} + L_{n-2}$, $n\ge 2$. We will use the following well-known facts about $F_n$, $L_n$, $\Gamma_n$, and $\Lambda_n$ without special mention:
\begin{enumerate}
\item $L_n\ =\ F_{n-1} + F_{n+1}$\ for $n \ge 1$,
\item $|V(\Gamma_n)|\ =\ F_{n+2}$\  for $n \ge 0$,
\item  $|E(\Gamma_n)|\ =\ (n F_{n+1} + 2(n+1)F_n)/5$\  for $n \ge 0$,
\item  $|V(\Lambda_n)|\ =\ L_{n}$\  for $n \ge 1$,\  $|V(\Lambda_0)|\ =\ 1$,
\item  $|E(\Lambda_n)|\ =\ n F_{n-1}$\  for $n \ge 0$.
\end{enumerate}
\begin{proposition}
\label{LucasBinomials}
\begin{eqnarray}
F_{n+1} &=& \sum_{k=0}^{\lfloor \frac{n}{2} \rfloor} \binom{n-k}{k}\ \ \mbox{for\ } n \ge -1,  \label{fbin} \\ 
L_n &=& \sum_{k=0}^{\lfloor \frac{n}{2} \rfloor} \frac{n}{n-k} \binom{n-k}{k}\ \,\mbox{for\ } n \ge 1,  \label{lbin} \\
\sum_{i=0}^{n} F_i L_{n-i} &=& (n+1)F_n \mbox{\ \ \, for\ } n \ge 0. \label{convol}
\end{eqnarray}
\end{proposition}
\proof Identities (\ref{fbin}) and (\ref{convol}) are well known. To prove (\ref{lbin}), note that by using (\ref{fbin}) twice and shifting the index of summation in the first sum,
\begin{eqnarray*}
L_n &=& F_{n-1} + F_{n+1}
\ =\  \sum_{k=0}^{\lfloor \frac{n-2}{2} \rfloor} \binom{n-k-2}{k} + \sum_{k=0}^{\lfloor \frac{n}{2} \rfloor} \binom{n-k}{k}\\
&=& \sum_{k=1}^{\lfloor \frac{n}{2} \rfloor} \binom{n-k-1}{k-1} + \sum_{k=0}^{\lfloor \frac{n}{2} \rfloor} \binom{n-k}{k}\\
&=& \sum_{k=0}^{\lfloor \frac{n}{2} \rfloor}\left(\frac{k}{n-k} + 1\right) \binom{n-k}{k}
\ =\ \sum_{k=0}^{\lfloor \frac{n}{2} \rfloor} \frac{n}{n-k} \binom{n-k}{k}. \qquad\mathqed
\end{eqnarray*}

\noindent
As usual, the automorphism group of a graph $G=(V,E)$ will be denoted by ${\rm Aut}(G)$.
The sets of orbits of ${\rm Aut}(G)$ acting on $V$ resp.\ $E$ will be denoted by $\mathcal{O}_V(G)$ resp.\ $\mathcal{O}_E(G)$, and their cardinalities by $o_V(G) = |\mathcal{O}_V(G)|$ and $o_E(G) = |\mathcal{O}_E(G)|$. In the latter case, we consider the natural action of ${\rm Aut}(G)$ on $E$, that is, for $g \in {\rm Aut}(G)$ the edge $\{u,v\}$ is mapped to $\{g(u),g(v)\}$. We denote the orbit of $u \in V$ resp.\ $e \in E$ under the action of ${\rm Aut}(G)$  by $\bar u$ resp.\ $\bar e$, and its cardinality by $|\bar u|$ resp.\ $|\bar e|$. In addition, we denote the number of orbits of size $k$ by
\begin{eqnarray*}
o_V(G,k) &=& |\{X \in \mathcal{O}_V(G);\ |X| = k\}|, \\
o_E(G,k) &=& |\{Y \in \mathcal{O}_E(G);\ |Y| = k\}|,
\end{eqnarray*}
so that the following identities hold:
\begin{eqnarray}
 \sum_k k\, o_V(G,k) = |V(G)|, && \sum_k o_V(G,k) = o_V(G), \label{abV} \\
 \sum_k k\, o_E(G,k) = |E(G)|, && \sum_k o_E(G,k) = o_E(G). \label{abE}
\end{eqnarray}
We denote the dihedral group of order $2n$ by $D_n$, and the set of fixed points of a group element $g$ acting on a set $A$, resp.\ its cardinality, by
\[
\mbox{Fix}_A(g)\ =\ \{u \in A;\ g(u) = u\}, \ \ \ \mbox{fix}_A(g)\ =\ |\mbox{Fix}_A(g)|.
\]
Finally, $\mathbb N$ is the set of all positive integers $\{1,2,3,\ldots\}$. For $n \in {\mathbb N}$,  $[n]$ denotes the set $\{1,2,\ldots, n\}$, and $[n]_0$ denotes the set $\{0,1,\ldots, n-1\}$.

%%%%%%%%%%%%%%%%%%%%%%%%%%%%%%%%%%%%%%%%%%%%%%%%%%%%%%%%%%%%%%%%%%%%%
%%%%%%%%%%%%%%%%%%%%%%%%%%%%%%%%%%%%%%%%%%%%%%%%%%%%%%%%%%%%%%%%%%%%%
\section{Dihedral transformations of strings}
\label{sec:strings}
%%%%%%%%%%%%%%%%%%%%%%%%%%%%%%%%%%%%%%%%%%%%%%%%%%%%%%%%%%%%%%%%%%%%%
%%%%%%%%%%%%%%%%%%%%%%%%%%%%%%%%%%%%%%%%%%%%%%%%%%%%%%%%%%%%%%%%%%%%%

As a preparation for what follows we investigate here some properties of dihedral transformations of nonempty strings defined over a finite alphabet. 

%%%%%%%%%%%%%%%%%%%%%%%%%%%%%%%%%%%%%%%%%%%%%%%%%%%%%%%%%%%%%%%%%%%%%
%\subsection{Primitive and asymmetric strings}
%\label{subsec:pastrings}
%%%%%%%%%%%%%%%%%%%%%%%%%%%%%%%%%%%%%%%%%%%%%%%%%%%%%%%%%%%%%%%%%%%%%

Let $\Sigma$ be an alphabet such that $0 \in \Sigma$. As usual, $\Sigma^n$ denotes the set of all strings of length $n$ over $\Sigma$ and $\Sigma^+ = \bigcup_{n=1}^\infty \Sigma^n$ denotes the set of all nonempty strings over $\Sigma$. If $u, v\in \Sigma^+$ and  $k \in {\mathbb N}$, we write $uv$ for the concatenation of $u$ and $v$, and $u^k$ for the concatenation of $k$ copies of $u$. For $u = u_1 u_2 \cdots u_n \in \Sigma^n$, we define its {\em length\/} $|u|$, {\em weight\/} $w(u) \in [n] \cup \{0\}$, {\em cyclic shift\/} $\alpha(u) \in \Sigma^n$, and {\em reversal\/} $\beta(u) \in \Sigma^n$ by
\begin{eqnarray*}
|u| &=& n, \\
w(u) &=& |\{i\in[n];\ u_i \ne 0\}|, \\
\alpha(u) &=& u_n u_1 \cdots u_{n-1}, \\
\beta(u) &=& u_n  u_{n-1} \cdots  u_1.
\end{eqnarray*}
It is straightforward to verify that
\begin{eqnarray}
x = \alpha^j(u) &\Longleftrightarrow& x_i = u_{(i-j) \bmod n} \ \Longleftrightarrow\ u_i = x_{(i+j) \bmod n}, \label{alphaj}\\
x = \beta(u) &\Longleftrightarrow& x_i = u_{(1-i) \bmod n} \ \Longleftrightarrow\ u_i = x_{(1-i) \bmod n},
\label{beta} \\
x = \alpha^j\beta(u) &\Longleftrightarrow& x_i = u_{(1-i+j) \bmod n} \ \Longleftrightarrow\ u_i = x_{(1-i+j) \bmod n}, 
\label{alphajbeta}
\end{eqnarray}
where $a \mbox{\ mod\ } n$ is the unique $i \in [n]$ such that $a \equiv i\ (\mbox{mod\ } n)$. Since $\alpha^n = \beta^2 = \mbox{id}$ and $\alpha\beta = \beta\alpha^{-1}$, the group generated by $\alpha$ and $\beta$ represents the action of $D_n$ on $\Sigma^n$. We denote the orbit of $u \in \Sigma^n$ under this action by $\bar u$, and its cardinality by $|\bar u|$.

\begin{lemma}
\label{ralpha}
For all $u \in \Sigma^+$, $j \in {\mathbb Z}$, and $k \in {\mathbb N}$, we have
\begin{enumerate}
\item[(i)] $\alpha^j(u^k) = (\alpha^j(u))^k$,
\item[(ii)] $\beta(u^k) = \beta(u)^k$.
\end{enumerate}
\end{lemma}

\proof This follows from (\ref{alphaj}), (\ref{beta}), and the fact that $(u^k)_i = u_{i \mbox{\ mod\ } |u|}$. \qed

%\begin{definition}
\noindent
For $u \in \Sigma^+$ define its {\em period} $p(u)$ and {\em exponent} $\ell(u)$ by
\begin{eqnarray*}
p(u) &=& \min\{k > 0;\ \alpha^k(u) = u\},\\
\ell(u) &=& \max\{k > 0;\ \exists v\in\Sigma^+\!: v^k = u\}.
\end{eqnarray*}
If $\ell(u) = 1$ then $u$ is {\em primitive}. It is well known \cite[Cor.\ 4.2]{lysc-1962} that for each $u$
there is a unique primitive string $\tau(u) \in \Sigma^+$ (called the {\em root} of $u$) such that $\tau(u)^{\ell(u)} =\, u$.
%\end{definition}

As an example, consider $\Sigma = \{0,1\}$ and the following strings $u \in \Sigma^4$:
\begin{itemize}
\item $u = 0000 = 0^4$: $p(u) = 1$, $\tau(u) = 0$, $\ell(u) = 4$, $u$ is not primitive,
\item $u = 0101 = (01)^2$: $p(u) = 2$, $\tau(u) = 01$, $\ell(u) = 2$, $u$ is not primitive,
\item $u = 0011 = (0011)^1$: $p(u) = 4$, $\tau(u) = 0011$, $\ell(u) = 1$, $u$ is primitive.
\end{itemize}

\begin{proposition}
\label{cu}
If $u \in \Sigma^n$, then 
\begin{enumerate}
\item[(i)] $p(u) = |\{u, \alpha(u), \alpha^2(u), \ldots, \alpha^{n-1}(u)\}|$,
\item[(ii)] $p(u^k) = p(u)$ for all $k \in {\mathbb N}$,
\item[(iii)] $p(u) \divides n$,
\item[(iv)] $p(u) = p(\tau(u)) = |\tau(u)|$\ and\ $\ell(u)\,p(u) = n$.
\end{enumerate}
\end{proposition}

\proof  
\begin{enumerate}
\item[(i)] Immediate from the definition of $p(u)$.
\item[(ii)] Since $u^k = v^k$ if and only if $u = v$, this follows from (i) and Lemma \ref{ralpha}(ii).
\item[(iii)] By (i), $p(u)$ is the size of the orbit of $u$ under the action of the cyclic group of order $n$ generated by $\alpha$,
hence $p(u)$ divides $n$.
\item[(iv)] Let $u = u_1 u_2 \cdots u_n$ and $v = u_1 u_2 \cdots u_{p(u)}$. From $\alpha^{p(u)}(u) = u$ it follows by (\ref{alphaj}) that $u_i = u_{i+p(u)}$ for $i = 1, 2, \ldots, n-p(u)$ (with indices taken mod $n$). By induction on $k$,  $u_i = u_{i+k p(u)}$ for $i = 1, 2, \ldots, n-k p(u)$ and $k = 0, 1, \ldots, \frac{n}{p(u)} - 1$. Hence $v = u_{1+kp(u)} u_{2+kp(u)} \cdots u_{p(u)+k p(u)}$ for $k = 0, 1, \ldots, \frac{n}{p(u)} - 1$, therefore $u = v ^{n/p(u)}$. If $u$ is primitive, this implies that $|u|/p(u) = 1$ and so $p(u) = |u|$. For arbitrary $u$ we then have, by (ii),
\[
p(u)\ =\ p(\tau(u)^{\ell(u)})\ =\ p(\tau(u))\ =\ |\tau(u)|,
\]
so $\ell(u)\,p(u) = \ell(u)|\tau(u)| = |\tau(u)^{\ell(u)}| = |u| = n$. \itemqed
\end{enumerate}

\begin{proposition}
\label{size2cu}
For each $u \in \Sigma^+$,
\[
|\bar u| = \left\{
\begin{array}{ll}
p(u), & \exists j:\ \beta(u) = \alpha^j(u), \\
2p(u), & \forall j:\ \beta(u) \ne \alpha^j(u).
\end{array}
\right.
\]
\end{proposition}

\proof Denote 
\begin{eqnarray*}
A(u) &=& \{u, \alpha(u), \alpha^2(u), \ldots, \alpha^{n-1}(u)\},\\
B(u) &=& \{\beta(u), \alpha \beta(u), \alpha^2 \beta(u), \ldots, \alpha^{n-1}\beta(u)\}.
\end{eqnarray*}
Then $\bar u = A(u) \cup B(u)$. If $\beta(u) = \alpha^j(u)$ for some $j$, it follows from $\alpha^{n} = \mbox{id}$ that $B(u) = A(u)$, hence $\bar u = A(u)$ and, by Proposition \ref{cu}(i), $|\bar u| = p(u)$.

Otherwise, if $\beta(u) \ne \alpha^j(u)$ for all $j$, it follows from $\alpha^{n} = \mbox{id}$ that $B(u) \cap A(u) = \emptyset$ and $|B(u)| = |A(u)|$, so, by Proposition \ref{cu}(i), $|\bar u| = 2\,|A(u)| = 2p(u)$. \qed

In order to determine the orbit size $|\bar u|$ more precisely, we distinguish between symmetric and asymmetric strings, where for $n \in {\mathbb N}$ a string $u \in \Sigma^n$ will be called:
\begin{itemize}
\item {\em symmetric\/} if $|\bar u| < 2n$, 
\item {\em asymmetric\/} if $|\bar u| = 2n$.
\end{itemize}

\begin{proposition}
\label{asyprim}
Every asymmetric string is primitive.
\end{proposition}

\proof Let $u \in \Sigma^n$ be asymmetric. Then $|\bar u| = 2n$, hence by Proposition \ref{size2cu}, $p(u) = n$. Then by Proposition \ref{cu}(iv), $\ell(u) = 1$ and $u$ is primitive. \qed

As an example, consider $\Sigma = \{0,1\}$ and the following strings $u \in \Sigma^6$:
\begin{itemize}
\item $u = 000000$ is symmetric and not primitive,
\item $u = 001100$ is symmetric and primitive,
\item $u = 010011$ is asymmetric (and hence primitive).
\end{itemize}

\begin{theorem}
\label{asym}
For each $u \in \Sigma^+$,
\[
|\bar u| = \left\{
\begin{array}{ll}
p(u), &  \tau(u) \mbox{\ symmetric}, \\
2p(u), & \tau(u) \mbox{\ asymmetric}.
\end{array}
\right.
\]
\end{theorem}

\proof Using Propositions \ref{size2cu}, \ref{cu}(iv) and Lemma \ref{ralpha} repeatedly we obtain
\begin{eqnarray*}
|\bar u| = 2p(u) &\Longleftrightarrow& \forall j:\,\beta(u) \ne \alpha^j(u)
\ \Longleftrightarrow\ \forall j:\,\beta(\tau(u)^{\ell(u)}) \ne \alpha^j(\tau(u)^{\ell(u)}) \\
&\Longleftrightarrow& \forall j:\,\beta(\tau(u))^{\ell(u)} \ne \alpha^j(\tau(u))^{\ell(u)} \\
&\Longleftrightarrow& \forall j:\,\beta(\tau(u)) \ne \alpha^j(\tau(u)) 
\ \Longleftrightarrow\ |\overline{\tau(u)}| = 2 p(\tau(u)) \\
&\Longleftrightarrow& |\overline{\tau(u)}| = 2 |\tau(u)|\ \Longleftrightarrow\ \tau(u) \mbox{\ asymmetric}. 
\end{eqnarray*} 
Together with Proposition \ref{size2cu} this proves the claim. \qed

\begin{corollary}
\label{ps}
A string $u \in \Sigma^n$ is primitive symmetric if and only if it is primitive and $\alpha^j \beta(u) = u$ for some $j \in [n]_0$.
\end{corollary}

\proof
By Theorem \ref{asym} and Proposition \ref{size2cu}, we have
\begin{eqnarray*}
u \mbox{\ primitive\ symmetric}
&\Longleftrightarrow& u \mbox{\ primitive} \,\land\, |\bar u| = p(u)\\
&\Longleftrightarrow& u \mbox{\ primitive} \,\land\, \exists j \in [n]_0:\ \beta(u) = \alpha^j(u)\\ 
&\Longleftrightarrow& u \mbox{\ primitive} \,\land\, \exists j \in [n]_0:\ \alpha^j \beta(u) = u. \qquad \mathqed
\end{eqnarray*}

%%%%%%%%%%%%%%%%%%%%%%%%%%%%%%%%%%%%%%%%%%%%%%%%%%%%%%%%%%%%%%%%%%%%%
%%%%%%%%%%%%%%%%%%%%%%%%%%%%%%%%%%%%%%%%%%%%%%%%%%%%%%%%
\section{Orbits of Fibonacci cubes}
\label{sec:orbit-fibonacci}
%%%%%%%%%%%%%%%%%%%%%%%%%%%%%%%%%%%%%%%%%%%%%%%%%%%%%%%%
%%%%%%%%%%%%%%%%%%%%%%%%%%%%%%%%%%%%%%%%%%%%%%%%%%%%%%%%%%%%%%%%%%%%%

According to \cite{cakl-2011}, for $n\ge 1$ the Fibonacci cube $\Gamma_n$ admits exactly one non-trivial automorphism, hence  $|{\rm Aut}(\Gamma_n)| = 2$ and the only orbit sizes are 1 and 2. For $n\ge 2$, the non-trivial automorphism coincides with the reversal map $\beta:V(\Gamma_n) \rightarrow V(\Gamma_n)$. We denote the set of Fibonacci strings of length $n$ which start with 0 resp.\ 1 by $V_0(\Gamma_n)$ resp.\  $V_1(\Gamma_n)$. We begin by enumerating Fibonacci palindromes.
\begin{proposition}
\label{fpal}
For $k \in {\mathbb N}$,
\[
\begin{array}{lllclll}
{\rm fix}_{V(\Gamma_{2k})}(\beta)&=&F_{k+1}, &\ & {\rm fix}_{V(\Gamma_{2k+1})}(\beta)&=&F_{k+3}, \\
{\rm fix}_{V_0(\Gamma_{2k})}(\beta)&=&F_k, && {\rm fix}_{V_0(\Gamma_{2k+1})}(\beta)&=&F_{k+2}, \\
{\rm fix}_{V_1(\Gamma_{2k})}(\beta)&=&F_{k-1}, && {\rm fix}_{V_1(\Gamma_{2k+1})}(\beta)&=&F_{k+1}.
\end{array}
\]
\end{proposition}

\proof
Let $u \in {\rm fix}_{V(\Gamma_{n})}(\beta)$. We distinguish four cases:
\begin{enumerate}
\item $n=2k$
\begin{enumerate}
\item $u$ starts with 0: $u = 0v00\beta(v)0$ with $v \in V(\Gamma_{k-2})$
\item $u$ starts with 1: $u = 10v00\beta(v)01$ with $v \in V(\Gamma_{k-3})$
\end{enumerate}
\item $n=2k+1$
\begin{enumerate}
\item $u$ starts with 0: $u = 0v0\beta(v)0$ with $v \in V(\Gamma_{k-1})$, or $u = 0v010\beta(v)0$ with $v \in V(\Gamma_{k-2})$
\item $u$ starts with 1: $u = 10v0\beta(v)01$ with $v \in V(\Gamma_{k-2})$, or $u = 10v010\beta(v)01$ with $v \in V(\Gamma_{k-3})$
\end{enumerate}
\end{enumerate}
The stated equalities now follow from $|V(\Gamma_k)| = F_{k+2}$ and $F_k + F_{k+1} = F_{k+2}$. \qed

%%%%%%%%%%%%%%%%%%%%%%%%%%%%%%%%%%%%%%%%%%%%%%%%%%%%%%%%
\subsection{Vertex orbits}
\label{sec:FibV}
%%%%%%%%%%%%%%%%%%%%%%%%%%%%%%%%%%%%%%%%%%%%%%%%%%%%%%%%

\begin{theorem}
\label{thm:number-of-orbits-Gamma}
Let $n \ge 2$. Then
\[
\begin{array}{lll}
o_V(\Gamma_n, 1) &=& F_{\lfloor \frac{n-(-1)^n}{2}\rfloor+2}, \\ 
o_V(\Gamma_n, 2) &=& \frac{1}{2} \left(F_{n+2} - F_{\lfloor \frac{n-(-1)^n}{2}\rfloor+2}\right), \\ 
o_V(\Gamma_n) &=& \frac{1}{2} \left(F_{n+2} + F_{\lfloor \frac{n-(-1)^n}{2}\rfloor+2}\right). 
\end{array}
\]
\end{theorem}

\proof
The orbits of size 1 correspond to the fixed points of $\beta$, i.e., to Fibonacci palindromes of length $n$, hence by Proposition \ref{fpal},
\[
o_V(\Gamma_n, 1) \ =\ \left\{
\begin{array}{ll}
F_{k+1}, & \mbox{if\ } n=2k, \\
F_{k+3}, & \mbox{if\ } n=2k+1,
\end{array}
\right.
\]
which can be combined into $F_{\lfloor \frac{n-(-1)^n}{2}\rfloor+2}$. The remaining two equalities now follow from (\ref{abV})  which in this case transforms into
\[
\begin{array}{lll}
o_V(\Gamma_n, 1) + 2 o_V(\Gamma_n, 2) &=& F_{n+2}, \\
o_V(\Gamma_n, 1) + o_V(\Gamma_n, 2) &=& o_V(\Gamma_n). \mathqed
\end{array}
\]
\begin{table}[h]
\[
\begin{array}{c|rrrrrrrrrrrrrrr}
n & 1& 2& 3& 4& 5& 6& 7& 8& 9& 10& 11& 12& 13& 14& 15 \\
\hline
|V(\Gamma_n)| & 2& 3& 5& 8& 13& 21& 34& 55& 89& 144& 233& 377& 610& 987& 1597 \\
o_V(\Gamma_n) & 1& 2& 4& 5& 9& 12& 21& 30& 51& 76& 127& 195& 322& 504& 826 \\
o_V(\Gamma_n, 1) & 0& 1& 3& 2& 5& 3& 8& 5& 13& 8& 21& 13& 34& 21& 55 \\
o_V(\Gamma_n, 2) & 1& 1& 1& 3& 4& 9& 13& 25& 38& 68& 106& 182& 288& 483& 771
\end{array}
\]
\caption{The numbers of vertices, all orbits, orbits of size 1, and orbits of size 2 in $V(\Gamma_n)$ for $n \le 15$ }
\end{table}

We remark that the numbers $o_V(\Gamma_{n})$ appear as solutions of other combinatorial enumeration problems as well. For instance, in~\cite{pago-1962}, the following problem is posed and solved: in how many ways can a $2\times (n+1)$ rectangle be tiled with dominoes (i.e., rectangles of sizes $2\times 1$ and $1\times 2$)? More precisely, the problem asks for the number of {\em distinct\/} tilings, where two tilings are considered distinct if one cannot be obtained from the other by reflections and rotations. As it turns out, the answer is given by $o_V(\Gamma_{n})$ (see also \cite[sequence A001224]{oeis}). For $n \in \{0,1\}$, this can be checked directly (note that the $2\times 2$  square has a single distinct tiling, due to the $90^{\circ}$ rotation). For $n \ge 2$, we present here a bijective proof of this fact:

Let $u \in V(\Gamma_n)$ and $v = u0 \in V(\Gamma_{n+1})$. Assign to $v$ a tiling of the $2 \times (n+1)$ rectangle with dominoes as follows. Going through $v$ from left to right, assign to each 0 a vertical domino and to each 10 a pair of horizontal dominoes. Conversely, to each tiling of the $2 \times (n+1)$ rectangle assign $v \in V(\Gamma_{n+1})$  by going through the tiling from left to right, coding vertical dominoes with 0 and pairs of horizontal dominoes with 10. Then $v$ ends with $0$; let $u\in V(\Gamma_n)$ be $v$ without the final 0. This establishes a bijection between $V(\Gamma_n)$ and the set of all tilings of the $2 \times (n+1)$ rectangle which preserves palindromes in both directions. Hence it gives rise to a bijection between ${\cal O}_V(\Gamma_n)$ and the set of all distinct tilings of the $2 \times (n+1)$ rectangle. 

Another family of combinatorial objects enumerated by $o_V(\Gamma_{n})$ are ordered integer partitions of $n+1$ with parts taken from the set $\{1, 2\}$ where two partitions are considered distinct if one cannot be obtained from the other by reflection. Such partitions are obviously in bijection with distinct domino tilings of the  $2\times (n+1)$ rectangle: to each part 1 in the partition assign a vertical domino, to each part 2 in the partition assign a pair of horizontal dominoes, and vice versa.

%%%%%%%%%%%%%%%%%%%%%%%%%%%%%%%%%%%%%%%%%%%%%%%%%%%%%%%%
\subsection{Edge orbits}
\label{sec:FibE}
%%%%%%%%%%%%%%%%%%%%%%%%%%%%%%%%%%%%%%%%%%%%%%%%%%%%%%%%

\begin{theorem} \label{T2}
For all $n\ge 0$,
\[
\begin{array}{lll}
o_E(\Gamma_n, 1) &=& \frac{1-(-1)^n}{2}  F_{\lfloor\frac{n+1}{2}\rfloor}, \\ 
o_E(\Gamma_n, 2) &=& \frac{1}{10} \left( nF_{n+1}+2(n+1)F_n\right) - \frac{1-(-1)^n}{4}  F_{\lfloor\frac{n+1}{2}\rfloor}, \\ 
o_E(\Gamma_n) &=& \frac{1}{10} \left( nF_{n+1}+2(n+1)F_n\right) + \frac{1-(-1)^n}{4}  F_{\lfloor\frac{n+1}{2}\rfloor}. 
\end{array}
\]
\end{theorem}

\proof
For $n \in \{0,1\}$, this can be checked directly. Let $n \ge 2$. The orbits of size 1 correspond to the fixed points of $\beta$, i.e., to
pairs $\{u,v\} \in E(\Gamma_n)$ such that $\{\beta(u),\beta(v)\} = \{u,v\}$. Since $u$ and $v$ differ in weight while $\beta$ preserves it, this is only possible if $\beta(u) = u$ and $\beta(v) = v$.  Hence $u$ and $v$ are Fibonacci palindromes of length $n$ differing in a single position. This is only possible if $n = 2k+1$ and
\begin{eqnarray*}
u &=& x000\beta(x), \\
v &=& x010\beta(x)
\end{eqnarray*}
or vice versa, for some $x\in V(\Gamma_{k-1})$. Hence
\begin{equation}
\label{fixeg}
o_E(\Gamma_n, 1)\ =\ \left\{
\begin{array}{ll}
|V(\Gamma_{k-1})| = F_{k+1}, & \mbox{if\ } n=2k+1, \\
0, & \mbox{if\ } n=2k,
\end{array}
\right.
\end{equation}
which can be written as $\frac{1-(-1)^n}{2}  F_{\lfloor\frac{n+1}{2}\rfloor}$. The remaining two equalities now follow from (\ref{abE})  which in this case transforms into
\[
\begin{array}{lll}
o_E(\Gamma_n, 1) + 2 o_E(\Gamma_n, 2) &=& \frac{1}{5} \left( nF_{n+1}+2(n+1)F_n\right), \\
o_E(\Gamma_n, 1) + o_E(\Gamma_n, 2) &=& o_E(\Gamma_n). \qquad \qquad \qquad \mathqed
\end{array}
\]
\begin{table}[h]
\[
\begin{array}{c|rrrrrrrrrrrrrr}
n & 1& 2& 3& 4& 5& 6& 7& 8& 9& 10& 11& 12& 13& 14 \\
\hline
|E(\Gamma_n)| & 1& 2& 5& 10& 20& 38& 71& 130& 235& 420& 744& 1308& 2285& 3970 \\
o_E(\Gamma_n) & 1& 1& 3& 5& 11& 19& 37& 65& 120& 210& 376& 654& 1149& 1985 \\
o_E(\Gamma_n, 1) &1& 0& 1& 0& 2& 0& 3& 0& 5& 0& 8& 0& 13& 0 \\
o_E(\Gamma_n, 2) & 0& 1& 2& 5& 9& 19& 34& 65& 115& 210& 368& 654& 1136& 1985
\end{array}
\]
\caption{The numbers of edges, all orbits, orbits of size 1, and orbits of size 2 in $E(\Gamma_n)$ for $n \le 14$ }
\end{table}

%%%%%%%%%%%%%%%%%%%%%%%%%%%%%%%%%%%%%%%%%%%%%%%%%%%%%%%%%%%%%%%%%%%%%
%%%%%%%%%%%%%%%%%%%%%%%%%%%%%%%%%%%%%%%%%%%%%%%%%%%%%%%%%%%%%%%%%%%%%
\section{Orbits of Lucas cubes}
\label{sec:Lucas}
%%%%%%%%%%%%%%%%%%%%%%%%%%%%%%%%%%%%%%%%%%%%%%%%%%%%%%%%%%%%%%%%%%%%%
%%%%%%%%%%%%%%%%%%%%%%%%%%%%%%%%%%%%%%%%%%%%%%%%%%%%%%%%%%%%%%%%%%%%%

Since Lucas cubes can be viewed as a symmetrization of Fibonacci cubes, the former should possess larger automorphism groups than the latter. 
Indeed, as shown in~\cite{cakl-2011}, ${\rm Aut}(\Lambda_n)$ is generated by the cyclic shift and reversal maps
$\alpha, \beta :V(\Lambda_n) \rightarrow V(\Lambda_n)$.
Hence ${\rm Aut}(\Lambda_n)$ = $\lbrace \mbox{id}, \alpha, \alpha^2, \cdots, \alpha^{n-1}, \beta, \alpha \beta, \alpha^2 \beta, \cdots,
\alpha^{n-1} \beta \rbrace \simeq D_n$ when  $n \ge 3$.

%%%%%%%%%%%%%%%%%%%%%%%%%%%%%%%%%%%%%%%%%%%%%%%%%%%%%%%%
\subsection{Vertex orbits}
\label{sec:LucasV}
%%%%%%%%%%%%%%%%%%%%%%%%%%%%%%%%%%%%%%%%%%%%%%%%%%%%%%%%

\begin{theorem} \label{thm:Lucas-vertex-orbits}
For all $n \in {\mathbb N}$,
\[
o_V(\Lambda_n)\ =\  
\frac{1}{2} \left(\frac{1}{n} \sum_{d | n} \varphi\left(\frac{n}{d}\right) L_{d} \ + \sum_{0 \le a \le \lfloor \frac{n}{2} \rfloor}\binom{\lfloor \frac{n}{2}\rfloor - \lfloor \frac{a+1}{2} \rfloor}{\lfloor \frac{a}{2} \rfloor}\right).
\]
\end{theorem}

\proof For $n \le 2$ this can be checked directly. Now let $n \ge 3$. Each $u \in V(\Lambda_n)$ can be constructed uniquely by concatenating $a$ copies of $01 \in V(\Lambda_2)$ and $n-2a$ copies of $0 \in V(\Lambda_1)$ for some $a$ between $0$ and $n/2$. Together this gives $n-a$ objects of two ``colors'' (01 and 0) arranged in a circle. Hence our problem is equivalent to determining the total number of 2-colorings of $n-a$ objects under the action of $D_{n-a}$, with $a$ objects of color 01 and $n-2a$ objects of color 0, for $0 \le a \le n/2$. From the well-known Redfield-P\'olya Theorem we have 
\begin{equation}
\label{polya}
o_V(\Lambda_n)\ =\ \sum_{a=0}^{\lfloor \frac{n}{2} \rfloor} [x_1^a x_2^{n-2a}]Z_{D_{n-a}}(x_1+x_2,\, x_1^2+x_2^2,\, \ldots,\, x_1^n+x_2^n)
\end{equation}
where 
\begin{equation}
\label{cindex}
Z_{D_n}(y_1,y_2,\ldots,y_n)\ =\ \frac{1}{2n}\sum_{d\divides n}\varphi(d)y_d^{n/d} + \frac{1}{4}
\left\{\begin{array}{ll}
2 y_1 y_2^{(n-1)/2}, & n \mbox{\ odd}, \\
y_1^2 y_2^{(n-2)/2} + y_2^{n/2}, & n \mbox{\ even},
\end{array}\right.
\end{equation}
is the cycle index polynomial of the natural action of $D_n$, and $[x_1^i x_2^j]p(x_1,x_2)$ denotes the coefficient of $x_1^i x_2^j$ in the polynomial $p(x_1,x_2)$.
From (\ref{polya}) and (\ref{cindex}) it follows by several applications of the Binomial Theorem that
\[
o_V(\Lambda_n)\ =\ \frac{1}{2} \left(c(n) +
\ \sum_{0 \le a \le \lfloor \frac{n}{2} \rfloor}\binom{\lfloor \frac{n}{2}\rfloor - \lfloor \frac{a+1}{2} \rfloor}{\lfloor \frac{a}{2} \rfloor}\right)
\]
where
\begin{equation}
\label{cn}
c(n)\ =\ \sum_{a=0}^{\lfloor \frac{n}{2} \rfloor} \frac{1}{n-a} \sum_{d | \gcd(n,a)} \varphi(d) \binom{\frac{n-a}{d}}{\frac{a}{d}}.
\end{equation}
By changing the order of summation on the right side of (\ref{cn}), writing $a = kd$, and replacing $d$ by $n/d$, we obtain
\begin{eqnarray*}
c(n) &=& \sum_{d | n} \varphi(d) \sum_{k=0}^{\lfloor \frac{n}{2d} \rfloor} \frac{1}{n-kd} \binom{\frac{n}{d}-k}{k} 
\ =\ \frac{1}{n} \sum_{d | n} \varphi\left(\frac{n}{d}\right) \sum_{k=0}^{\lfloor \frac{d}{2} \rfloor} \frac{d}{d-k} \binom{d-k}{k}.
\end{eqnarray*}
From Proposition \ref{LucasBinomials} it now follows that 
\[
c(n)\ =\ \frac{1}{n} \sum_{d | n} \varphi\left(\frac{n}{d}\right)L_d,
\]
proving the claim. 
\qed

\noindent
We remark that there is a natural bijection between the set of Lucas strings of length $n$ and the set of circular strings of $n$ black and white beads with no two black beads adjacent. Hence $o_V(\Lambda_n)$ counts black-and-white bracelets with no consecutive black beads. These numbers are known (cf.\ \cite[sequence A129526]{oeis}), as are the numbers $c(n)$ which count  black-and-white necklaces with no consecutive black beads  (cf.\ \cite[sequence A000358]{oeis}).

Now we determine the sizes of orbits. They must divide $2n$, but the following result shows that not all divisors of $2n$ appear as orbit sizes. 

\begin{theorem}
\label{LucasVorbitsize}
For all $n \in {\mathbb N}$,
\[
\{|X|;\ X \in {\cal O}_V(\Lambda_n)\}\ =\ 
\{k \ge 1;\ k \divides n\} \cup\, \{k \ge 18;\ k \divides 2n\}.
\]
\end{theorem}
To prove the theorem, we first show: 

\begin{lemma}
\label{length9}
$V(\Lambda_n)$ contains an asymmetric string if and only if $n \ge 9$.
\end{lemma}

\proof To prove that all $u \in V(\Lambda_n)$ with $n<9$ are symmetric we order strings according to their weight $w(u)$. Since conjugate strings are either all symmetric or all asymmetric, it suffices to consider a single representative from each orbit:
\begin{itemize}
\item $w(u) = 0$: $u = 0^n$ is symmetric since $\beta(u) = u$.
\item $w(u) = 1$: $u = 10^{n-1}$ is symmetric since $\beta(u) = 0^{n-1}1 = \alpha^{-1}(u)$.
\item $w(u) = 2$: $u = 10^a 10^b$ where $a, b \ge 1$ is symmetric since $\beta(u) = 0^b 10^a 1 = \alpha^b(u)$.
\item $w(u) = 3$: $u = 10^a 10^b 10^c$ where $a, b, c \ge 1$. If any two of $a,b,c$ are equal, then $u$ is symmetric since
\begin{itemize}
\item if $a = b$: $u = 10^a 10^a 10^c$, $\beta(u) = 0^c 10^a 10^a 1 = \alpha^c(u)$,
\item if $a = c$: $u = 10^a 10^b 10^a$, $\beta(u) = 0^a 10^b 10^a 1 = \alpha^{-1}(u)$,
\item if $b = c$: $u = 10^a 10^b 10^b$, $\beta(u) = 0^b 10^b 10^a 1 = \alpha^{-a-2}(u)$.
\end{itemize}
Hence in any asymmetric $u$ of this form we have $a+b+c \ge 1+2+3 = 6$ and $n = a+b+c+3 \ge 9$. 
\item $w(u) = 4$: $u = 10^a 10^b 10^c 10^d$ where $a, b, c, d \ge 1$. If $a=b=c=d=1$ then $u = 10101010$ is symmetric since $\beta(u) = \alpha(u)$.
Otherwise $a+b+c+d \ge 5$ and $n = a+b+c+d+4 \ge 9$.
\item $w(u)\ge 5$: any $u \in V(\Lambda_n)$ with $w(u)\ge 5$ contains at least 5 zeros, so in this case $n \ge w(u) + 5 \ge 10$.
\end{itemize}

Conversely, assume that $n \ge 9$ and let $u = 1010010^{n-6}$. Then $u \in V(\Lambda_n)$, $u$ is obviously primitive, and $\beta(u) = 0^{n-6}100101$.
Since $n-6 \ge 3$, the only $j \in [n-1]$ such that $\alpha^j(u)$ starts with $0^{n-6}$ is $n-6$, but $\alpha^{n-6}(u) = 0^{n-6}101001
\ne \beta(u)$. Hence $u$ is asymmetric. \qed

\noindent{\bf Proof of Theorem~\ref{LucasVorbitsize}.\/}
% Clearly, $\{|\bar u|;\ u \in V(\Lambda_n)\}\ \subseteq\ \{k \ge 1;\ k \divides 2n\}$.
For $n \in \{0,1,2\}$, this can be checked directly. Let $n \ge 3$. Which divisors $k$ of $2n$ appear as orbit sizes in the action of ${\rm Aut}(\Lambda_n)$ on $V(\Lambda_n)$? We distinguish two cases:

{\em Case 1:} $k \divides n$. If $k = 1$, take $u = 0^n$; then $u \in V(\Lambda_n)$ and $|\bar u| = 1$. For $k \ge 2$, take $u = (10^{k-1})^{n/k} \in V(\Lambda_n)$. Since $10^{k-1}$ is primitive, $\ell(u) = n/k$ and $p(u) = k$, by Proposition \ref{cu}(iv). Since $\beta(u) = (0^{k-1}1)^{n/k} = \alpha^{k-1}(u)$, it follows from Proposition \ref{size2cu} that $|\bar u| = p(u) = k$. This shows that $\{k \ge 1;\ k \divides n\} \subseteq \{|\bar u|;\ u \in V(\Lambda_n)\}$.

{\em Case 2:\ $k \divides 2n$, but $k \notdiv n$}. Here $k$ is even. Note that it suffices to prove that there is an orbit of size $k$ 
in $V(\Lambda_n)$ if and only if $k \ge 18$.

Assume first that $u \in V(\Lambda_n)$ is such that $|\bar u| = k$. Then $k \in \{p(u), 2p(u)\}$ by Proposition \ref{size2cu}. 
Since $p(u) \divides n$ by Proposition \ref{cu}(iii) but $k \notdiv n$, we conclude that $k = 2p(u)$.
Theorem \ref{asym} now implies that $\tau(u)$ is asymmetric. By Lemma \ref{length9}, $|\tau(u)| \ge 9$, hence by Proposition \ref{cu}(iv), $p(u) \ge 9$ as well, so $k \ge 18$.

Conversely, assume that $k \ge 18$. By Lemma \ref{length9}, there is an asymmetric $v \in V(\Lambda_{k/2})$. Let $u = v^s$
where $s = 2n/k$. Then $u \in V(\Lambda_n)$ and by Proposition \ref{asyprim}, $v = \tau(u)$. By Theorem \ref{asym} and Proposition \ref{cu}(iv), $|\bar u| = 2p(u) = 2|v| = k$. \qed

In order to determine the number of orbits of size $k$, we enumerate primitive symmetric and asymmetric Lucas strings. We denote the set of primitive Lucas strings of length $n$ by $V_p(\Lambda_n)$. For $n \in {\mathbb N}$ we define
\[
\begin{array}{lll}
p_n &=& |V_p(\Lambda_n)| \ =\ |\{u \in V(\Lambda_n);\ u \mbox{\ primitive}\}|, \\
s_n &=& |\{u \in V(\Lambda_n);\ u \mbox{\ primitive\ symmetric}\}|, \\
t_n &=& \sum_{j=0}^{n-1}\mbox{fix}_{V_p(\Lambda_n)}(\alpha^j \beta), \\
a_n &=& |\{u \in V(\Lambda_n);\ u \mbox{\ asymmetric}\}.
\end{array}
\]

\begin{lemma}
\label{uniquefp}
Let $u\in \Sigma^n$ and $\alpha^j\beta(u) = \alpha^k\beta(u)$ for some $j, k \in [n]_0$ with $j<k$. Then $u$ is not primitive.
\end{lemma}

\proof Since $k-j > 0$, we have by Lemma \ref{ralpha}(ii),
\[
\begin{array}{lllll}
\alpha^j\beta(u) = \alpha^k\beta(u) &\Longrightarrow& \alpha^{k-j}\beta(u) = \beta(u) &\Longrightarrow& p(\beta(u)) \le k-j < n \\
 &\Longrightarrow& \beta(u) \mbox{\ not\ primitive\ } &\Longrightarrow& u \mbox{\ not\ primitive}. \mathqed
\end{array}
\]

\begin{corollary}
\label{ts}
\begin{itemize}
\item[(i)] The sets $\mbox{Fix}_{V_p(\Lambda_n)}(\alpha^j\beta)$, for $j \in [n]_0$, are pairwise disjoint.
\item[(ii)] For all $n \in {\mathbb N}$, $s_n = t_n$.
\end{itemize}
\end{corollary}

\proof
 \begin{itemize}
\item[(i)] Assume that $u \in \mbox{Fix}_{V_p(\Lambda_n)}(\alpha^j\beta) \cap\mbox{Fix}_{V_p(\Lambda_n)}(\alpha^k\beta)$. Then $\alpha^j\beta(u) = u = \alpha^k\beta(u)$ and $u$ is primitive, hence Lemma \ref{uniquefp} implies that $j = k$, proving (i). 
\item[(ii)] Using (i) and Corollary \ref{ps}, we obtain
\[
t_n = \sum_{j=0}^{n-1}|\mbox{Fix}_{V_p(\Lambda_n)}(\alpha^j \beta)| = |\bigcup_{j=0}^{n-1}\mbox{Fix}_{V_p(\Lambda_n)}(\alpha^j \beta)| = s_n.
\qquad\quad\mathqed
\]
\end{itemize}

\begin{lemma}
\label{palindr}
Let $u\in \Sigma^n$ and $0 \le j < n$. Then $\alpha^j\beta(u) =u$ if and only if there are $x \in \Sigma^j$ and $y \in \Sigma^{n-j}$ such that $x = \beta(x)$, $y = \beta(y)$, and $xy = u$.
\end{lemma}

\proof Assuming that $\alpha^j\beta(u) = u$, let $x = u_1 \cdots u_j$ and $y = u_{j+1} \cdots u_n$. Then $xy = u = \alpha^j\beta(u) = \alpha^j\beta(xy) = \alpha^j(\beta(y)\beta(x)) = \beta(x)\beta(y)$, hence $\beta(x) = x$ and $\beta(y) = y$. 

Conversely, assuming that  $x \in \Sigma^j$ and $y \in \Sigma^{n-j}$ are such that  $x = \beta(x)$, $y = \beta(y)$, and $xy = u$, we have  $\alpha^j\beta(u) =  \alpha^j\beta(xy) = \alpha^j(\beta(y)\beta(x)) = \beta(x)\beta(y) = xy = u$. \qed

In the next theorem, we express the numbers of primitive, primitive symmetric, and asymmetric Lucas strings of length $n$ by means of M\"obius function.

\begin{theorem}
\label{psa}
For all $n \in {\mathbb N}$,
\begin{itemize}
\item[(i)] $p_n\ =\ \displaystyle\sum_{d\divides n}\mu\left(\frac{n}{d}\right) L_d$,
\item[(ii)] $s_n\ =\ n\displaystyle\sum_{d\divides n}\mu\left(\frac{n}{d}\right) F_{\lfloor \frac{d}{2}\rfloor + 2}$,
\item[(iii)] $a_n\ =\ \displaystyle\sum_{d\divides n}\mu\left(\frac{n}{d}\right) (L_d - n  F_{\lfloor \frac{d}{2}\rfloor + 2})$.
\end{itemize}
\end{theorem}

\proof
If $u$ is a Lucas string and $k\in {\mathbb N}$, then $\tau(u)$ and $u^k$ are Lucas strings as well. Hence we can enumerate $u \in V(\Lambda_n)$ by $|\tau(u)|$ which is a divisor of $n$. This yields
\[
L_n = |V(\Lambda_n)| = \sum_{d\divides n} p_d, 
\]
from which (i) follows by M\"obius inversion.

To derive (ii), we compute $t_n$ and invoke Corollary \ref{ts}(ii). By Lemma \ref{ralpha}, 
\begin{eqnarray*}
\mbox{fix}_{V(\Lambda_n)}(\alpha^j \beta) &=& |\{u \in V(\Lambda_n);\ \alpha^j\beta(u) = u\}| \\ 
&=& |\bigcup_{d\divides n} \{u \in V(\Lambda_n);\ |\tau(u)| = d \,\land\, \alpha^j\beta(\tau(u)) = \tau(u)\}| \\ 
&=& |\bigcup_{d\divides n} \{v \in V(\Lambda_d);\ v \mbox{\ primitive} \,\land\, \alpha^j\beta(v) = v\}| \\ 
&=& |\bigcup_{d\divides n} \mbox{Fix}_{V_p(\Lambda_d)}(\alpha^j \beta)| \ =\  \sum_{d\divides n} \mbox{fix}_{V_p(\Lambda_d)}(\alpha^j \beta),
\end{eqnarray*}
hence by M\"obius inversion,
\[
 \mbox{fix}_{V_p(\Lambda_n)}(\alpha^j \beta)\ =\ \sum_{d\divides n} \mu\left(\frac{n}{d}\right)  \mbox{fix}_{V(\Lambda_d)}(\alpha^j \beta),
\]
and, by definition of $t_n$,
\begin{eqnarray}
t_n &=& \sum_{j=0}^{n-1} \sum_{d\divides n} \mu\left(\frac{n}{d}\right)  \mbox{fix}_{V(\Lambda_d)}(\alpha^j \beta) 
\ =\  \sum_{d\divides n} \mu\left(\frac{n}{d}\right) \sum_{j=0}^{n-1} \mbox{fix}_{V(\Lambda_d)}(\alpha^j \beta) \nonumber \\
 &=&  \sum_{d\divides n} \mu\left(\frac{n}{d}\right) \frac{n}{d} \sum_{j=0}^{d-1} \mbox{fix}_{V(\Lambda_d)}(\alpha^j \beta). \label{tn}
\end{eqnarray}
To evaluate $\sum_{j=0}^{d-1} \mbox{fix}_{V(\Lambda_d)}(\alpha^j \beta)$, we use Lemma \ref{palindr} according to which we need to count strings $xy \in V(\Lambda_d)$ where $x \in \mbox{Fix}_{V(\Gamma_j)}(\beta)$ and $y \in \mbox{Fix}_{V(\Gamma_{d-j})}(\beta)$ are Fibonacci palindromes. We distinguish four cases according to the parities of $d$ and $j$, computing the numbers of such strings in each case by means of Proposition \ref{fpal}:
\begin{enumerate}
\item $d = 2k$
\begin{enumerate}
\item $j = 2i$
\begin{itemize}
\item $x  \in  \mbox{Fix}_{V_0(\Gamma_{2i})}(\beta)$, $y \in \mbox{Fix}_{V_0(\Gamma_{2(k-i)})}(\beta)$:\ \ $F_i F_{k-i}$ strings
\item $x  \in  \mbox{Fix}_{V_0(\Gamma_{2i})}(\beta)$, $y \in \mbox{Fix}_{V_1(\Gamma_{2(k-i)})}(\beta)$:\ \ $F_i F_{k-i-1}$ strings
\item $x  \in  \mbox{Fix}_{V_1(\Gamma_{2i})}(\beta)$, $y \in \mbox{Fix}_{V_0(\Gamma_{2(k-i)})}(\beta)$:\ \ $F_{i-1} F_{k-i}$ strings
\end{itemize}
in all: $F_i F_{k-i+1} + F_{i-1} F_{k-i} = \frac{1}{2}(F_i L_{k-i} + L_i F_{k-i})$ strings

\item $j = 2i+1$
\begin{itemize}
\item $x  \in  \mbox{Fix}_{V_0(\Gamma_{2i+1})}(\beta)$, $y \in \mbox{Fix}_{V_0(\Gamma_{2(k-i-1)+1})}(\beta)$:\ \ $F_{i+2} F_{k-i+1}$ strings
\item $x  \in  \mbox{Fix}_{V_0(\Gamma_{2i+1})}(\beta)$, $y \in \mbox{Fix}_{V_1(\Gamma_{2(k-i-1)+1})}(\beta)$:\ \ $F_{i+2} F_{k-i}$ strings
\item $x  \in  \mbox{Fix}_{V_1(\Gamma_{2i+1})}(\beta)$, $y \in \mbox{Fix}_{V_0(\Gamma_{2(k-i-1)+1})}(\beta)$:\ \ $F_{i+1} F_{k-i+1}$ strings
\end{itemize}
in all: $F_{i+2} F_{k-i+2} + F_{i+1} F_{k-i+1} = \frac{1}{2}(F_{i+2} L_{k-i+1} + L_{i+2} F_{k-i+1})$ strings
\end{enumerate}

\item $d = 2k+1$
\begin{enumerate}
\item $j = 2i$
\begin{itemize}
\item $x  \in  \mbox{Fix}_{V_0(\Gamma_{2i})}(\beta)$, $y \in \mbox{Fix}_{V_0(\Gamma_{2(k-i)+1})}(\beta)$:\ \ $F_i F_{k-i+2}$ strings
\item $x  \in  \mbox{Fix}_{V_0(\Gamma_{2i})}(\beta)$, $y \in \mbox{Fix}_{V_1(\Gamma_{2(k-i)+1})}(\beta)$:\ \ $F_i F_{k-i+1}$ strings
\item $x  \in  \mbox{Fix}_{V_1(\Gamma_{2i})}(\beta)$, $y \in \mbox{Fix}_{V_0(\Gamma_{2(k-i)+1})}(\beta)$:\ \ $F_{i-1} F_{k-i+2}$ strings
\end{itemize}
in all: $F_i F_{k-i+3} + F_{i-1} F_{k-i+2} = \frac{1}{2}(F_i L_{k-i+2} + L_i F_{k-i+2})$ strings

\item $j = 2i+1$
\begin{itemize}
\item $x  \in  \mbox{Fix}_{V_0(\Gamma_{2i+1})}(\beta)$, $y \in \mbox{Fix}_{V_0(\Gamma_{2(k-i)})}(\beta)$:\ \ $F_{i+2} F_{k-i}$ strings
\item $x  \in  \mbox{Fix}_{V_0(\Gamma_{2i+1})}(\beta)$, $y \in \mbox{Fix}_{V_1(\Gamma_{2(k-i)})}(\beta)$:\ \ $F_{i+2} F_{k-i-1}$ strings
\item $x  \in  \mbox{Fix}_{V_1(\Gamma_{2i+1})}(\beta)$, $y \in \mbox{Fix}_{V_0(\Gamma_{2(k-i)})}(\beta)$:\ \ $F_{i+1} F_{k-i}$ strings
\end{itemize}
in all: $F_{i+2} F_{k-i+1} + F_{i+1} F_{k-i} = \frac{1}{2}(F_{i+2} L_{k-i} + L_{i+2} F_{k-i})$ strings

\end{enumerate}
\end{enumerate}
Again we distinguish two cases according to the parity of $d$, and split the sum into two according to the parity of the summation index:

\begin{enumerate}
\item $d=2k$:
\begin{eqnarray*}
\lefteqn{\sum_{j=0}^{d-1} \mbox{fix}_{V(\Lambda_d)}(\alpha^j \beta) 
\ =\ \sum_{i=0}^{k-1} \mbox{fix}_{V(\Lambda_d)}(\alpha^{2i} \beta) + \sum_{i=0}^{k-1} \mbox{fix}_{V(\Lambda_d)}(\alpha^{2i+1} \beta)} \\
&=& \frac{1}{2} \left( \sum_{i=0}^{k-1} F_i L_{k-i} + \sum_{i=0}^{k-1} L_i F_{k-i} +  \sum_{i=0}^{k-1}F_{i+2} L_{k-i+1} +  \sum_{i=0}^{k-1}L_{i+2} F_{k-i+1}\right) \\
&=& 2k F_{k+2}\ =\ d F_{\lfloor\frac{d}{2}\rfloor + 2},
\end{eqnarray*}
by shifting summation indices and applying (\ref{convol}) to each of the four sums.

\item $d=2k+1$:
\begin{eqnarray*}
\lefteqn{\sum_{j=0}^{d-1} \mbox{fix}_{V(\Lambda_d)}(\alpha^j \beta) 
\ =\ \sum_{i=0}^k \mbox{fix}_{V(\Lambda_d)}(\alpha^{2i} \beta) + \sum_{i=0}^{k-1} \mbox{fix}_{V(\Lambda_d)}(\alpha^{2i+1} \beta)} \\
&=& \frac{1}{2} \left( \sum_{i=0}^k F_i L_{k-i+2} +\sum_{i=0}^k L_i F_{k-i+2} +\sum_{i=0}^{k-1} F_{i+2} L_{k-i} +\sum_{i=0}^{k-1} L_{i+2} F_{k-i} \right) \\
&=& (2k+1) F_{k+2}\ =\ d F_{\lfloor\frac{d}{2}\rfloor + 2},
\end{eqnarray*}
by shifting summation indices and applying (\ref{convol}) to each of the four sums.
\end{enumerate}
The final expression is the same in both cases, so by (\ref{tn}) we obtain
\[
s_n\ =\ t_n\ =\ \sum_{d\divides n} \mu\left(\frac{n}{d}\right) \frac{n}{d} \sum_{j=0}^{d-1} \mbox{fix}_{V(\Lambda_d)}(\alpha^j \beta)\ =\ n \sum_{d\divides n} \mu\left(\frac{n}{d}\right)  F_{\lfloor\frac{d}{2}\rfloor + 2},
\]
proving (ii).

Finally, (iii) follows from (i) and (ii) by noting that $a_n = p_n - s_n$. \qed
\begin{table}[h]
\[
\begin{array}{c|rrrrrrrrrrrrrrrr}
n & 1& 2& 3& 4& 5& 6& 7& 8& 9& 10& 11& 12& 13& 14& 15& 16 \\
\hline
L_n & 1& 3& 4& 7& 11& 18& 29& 47& 76& 123& 199& 322& 521& 843& 1364& 2207 \\
p_n & 1& 2& 3& 4& 10& 12& 28& 40& 72& 110& 198& 300& 520& 812& 1350& 2160 \\
s_n & 1& 2& 3& 4& 10& 12& 28& 40& 54& 90& 132& 180& 260& 392& 450& 752 \\
a_n & 0& 0& 0& 0& 0& 0& 0& 0& 18& 20& 66& 120& 260& 420& 900& 1408
\end{array}
\]
\caption{The numbers of all  Lucas strings, primitive Lucas strings, symmetric primitive Lucas strings, and asymmetric Lucas strings of length $n \le 16$ }
\end{table}

\begin{theorem}
\label{orbitsizes}
For all $n \in {\mathbb N}$ and $k \divides 2n$,
\[
o_V(\Lambda_n, k)\ =\ \left\{
\begin{array}{ll}
\displaystyle\sum_{d\divides k}\mu\left(\frac{k}{d}\right) F_{\lfloor \frac{d}{2}\rfloor + 2}, & k \divides n \ \land\, k {\rm\ odd} \\
\displaystyle\sum_{d\divides k}\mu\left(\frac{k}{d}\right) F_{\lfloor \frac{d}{2}\rfloor + 2} + \frac{1}{k} \displaystyle\sum_{d\divides \frac{k}{2}}\mu\left(\frac{k}{2d}\right) (L_d - \frac{k}{2}  F_{\lfloor \frac{d}{2}\rfloor + 2}), & k \divides n \ \land\, k {\rm\ even} \\
\displaystyle\frac{1}{k} \sum_{d\divides \frac{k}{2}}\mu\left(\frac{k}{2d}\right) (L_d - \frac{k}{2}  F_{\lfloor \frac{d}{2}\rfloor + 2}), & k \divides 2n \,\land\, k \notdiv n
\end{array}
\right.
\]
\end{theorem}

\proof
By Theorem \ref{asym}, 
\[
|\bar u| = k\ \Longleftrightarrow\ (|\tau(u)| = k \land \tau(u) \mbox{\ symmetric}) \lor  (2|\tau(u)| = k \land \tau(u) \mbox{\ asymmetric}).
\]
Since $\tau(u)$ is primitive and $|\tau(u)|$ divides $n$, it follows that
\begin{eqnarray*}
|\bar u| = k &\Longleftrightarrow& \left(k \divides n\ \land\ \exists v \in V(\Lambda_k)\!\!: (v \mbox{\ symmetric\ primitive}\,\land\, u = v^{n/k})\right) \\
&& \lor\ \left(k \divides 2n\ \land\ k \mbox{\ even}\ \land\ \exists v \in V(\Lambda_{k/2})\!\!: (v \mbox{\ asymmetric} \,\land\, u = v^{2n/k})\right),
\end{eqnarray*}
hence
\begin{equation}
\label{sak}
o_V(\Lambda_n, k)\ =\ \frac{1}{k} \left\{
\begin{array}{ll}
s_k, & k \divides n \ \land\, k {\rm\ odd}, \\
s_k + a_{k/2}, & k \divides n \ \land\, k {\rm\ even}, \\
a_{k/2}, & k \divides 2n \,\land\, k \notdiv n.
\end{array}
\right.
\end{equation}
Together with Proposition \ref{psa}, this yields the statement of the theorem. \qed

At first sight, it seems surprising that the summation formulas in Theorem \ref{orbitsizes} are free of $n$, therefore we try to explain this phenomenon here. Since both $\alpha$ and $\beta$ commute with taking powers,  the size of the orbit of $u$ equals the size of the orbit of $u^k$, which means that $|\bar u|$ depends only on $|\tau(u)|$. Primitive strings of length $n$ are of two types: symmetric (with orbit size $n$) and asymmetric (of orbit size $2n$). Hence in order to enumerate orbits of size $k$ in $V(\Lambda_n)$, we need to enumerate primitive symmetric Lucas strings of length $k$ and asymmetric Lucas strings of length $k/2$, with the latter possibility applicable only to even $k$. Furthermore, we wish to obtain a string of length $n$ as a power of our primitive string, therefore $n/k$ (in the first case) resp.\ $2n/k$ (in the second case) must be an integer. In summary, if $k$ is even, only the first case applies, hence we obtain $s_k$ such strings. If $k$ is even and divides $n$, both cases apply and there are $s_k + a_{k/2}$ such strings. If $k$ is even but does not divide $n$  (however it must divide $2n$), only the second case applies, yielding $a_{k/2}$ strings. Finally, since we count orbits rather than individual strings, we divide by the orbit size $k$ and obtain \eqref{sak}.
\begin{table}[h]
\[
\begin{array}{c|rrrrrrrrrrrrrrrrrr}
n & 1& 2& 3& 4& 5& 6& 7& 8& 9& 10& 11& 12& 13& 14& 15& 16& 17 &18 \\
\hline
o_V(\Lambda_n) & 1 &  2 &  2 &  3 &  3 &  5 &  5 &  8 &  9 &  14 &  16 &  26 &  31 &  49 &  64 &  99 &  133 &  209 \\
o_V(\Lambda_n, n) & 1 &  1 &  1 &  1 &  2 &  2 &  4 &  5 &  6 &  9 &  12 &  15 &  20 &  28 &  30 &  47 &  54 &  79 \\
o_V(\Lambda_n, 2n) & 0 &  0 &  0 &  0 &  0 &  0 &  0 &  0 &  1 &  1 &  3 &  5 &  10 &  15 &  30 &  44 &  78 &  119
\end{array}
\]
\caption{The numbers of all orbits, orbits of size $n$, and orbits of size $2n$ in $V(\Lambda_n)$ for $n \le18$ }
\end{table}

%%%%%%%%%%%%%%%%%%%%%%%%%%%%%%%%%%%%%%%%%%%%%%%%%%%%%%%%
\subsection{Edge orbits}
%%%%%%%%%%%%%%%%%%%%%%%%%%%%%%%%%%%%%%%%%%%%%%%%%%%%%%%%

We now present a surprising relationship between $o_E(\Lambda_n)$ and $o_V(\Gamma_{n}),$ namely that the former sequence is just a shift of the latter. 

\begin{theorem}
\begin{itemize}
\item[(i)] For $n \ge 5$, $o_E(\Lambda_n)\ =\ o_V(\Gamma_{n-3})$.
\item[(ii)] For all $n \in {\mathbb N}$,
\[
o_E(\Lambda_n)\ =\ \frac{1}{2} \left(F_{n-1} + F_{\lfloor \frac{n+1+(-1)^n}{2}\rfloor}\right).
\]
\end{itemize}
\end{theorem}

\proof
To prove (i), assume that $n \ge 5$. Then ${\rm Aut}(\Gamma_{n-3}) = \{\mbox{id}, \beta\}$ and ${\rm Aut}(\Lambda_n)$ = $\lbrace \mbox{id}, \alpha, \cdots, \alpha^{n-1}$, $\beta, \alpha \beta, \cdots, \alpha^{n-1} \beta \rbrace$ $\simeq$ $D_n$. Define $s\!:\ E(\Lambda_n) \to V(\Gamma_{n-3})$ as follows (all indices will be taken mod $n$): For $e = \{u,v\} \in E(\Lambda_n)$, let $i \in [n]$ be such that $u_i \ne v_i$ where $u = u_1 u_2 \cdots u_n$, $v = v_1 v_2 \ldots v_n$. Then $\{u_i,v_i\} = \{0,1\}$ and $u_{i-1} = v_{i-1} = u_{i+1} = v_{i+1} = 0$. Now set
\[
s(e)\ =\ u_{i+2}u_{i+3} \cdots u_{i+n-2}.
\]
Note that $s(e)$ may contain $u_n u_1$ as a substring, but as $u$ is a Lucas string,  $u_n u_1 \ne 11$, hence $s(e) \in V(\Gamma_{n-3})$ as desired.  Let  $f = \{x,y\} \in E(\Lambda_n)$ be such that $\bar e = \bar f$. We claim that $\overline{s(e)} = \overline{s(f)}$. To prove this, we distinguish two cases:

{\em Case 1:} $f = \alpha^j(e) = \{\alpha^j(u), \alpha^j(v)\}$ for some $j$. Without loss of generality assume that $x = \alpha^j(u)$ and $y = \alpha^j(v)$. By (\ref{alphaj}), $x_{i+j} \ne y_{i+j}$, hence
\[
s(f) = x_{i+j+2}x_{i+j+3} \cdots x_{i+j+n-2} = \ u_{i+2}u_{i+3} \cdots u_{i+n-2} = s(e).
\]

{\em Case 2:} $f = \alpha^j\beta(e) = \{\alpha^j\beta(u), \alpha^j\beta(v)\}$ for some $j$. Without loss of generality assume that $x = \alpha^j\beta(u)$ and $y = \alpha^j\beta(v)$. By (\ref{alphajbeta}), $x_{1-i+j} \ne y_{1-i+j}$, hence
\begin{eqnarray*}
s(f) &=& x_{3-i+j}x_{4-i+j} \cdots x_{n-1-i+j}\ =\ u_{i-2}u_{i-3} \cdots u_{i-n+2} \\
       &=& \beta(u_{i-n+2}\cdots  u_{i-3}u_{i-2})\ =\ \beta(u_{i+2}\cdots  u_{i+n-3}u_{i+n-2})\ =\ \beta(s(e)).
\end{eqnarray*}
In either case, $\overline{s(e)} = \overline{s(f)}$ as claimed. This implies that $s$ induces a well-defined mapping $\tilde s\!:\ \mathcal{O}_E(\Lambda_n) \to \mathcal{O}_V(\Gamma_{n-3})$ which satisfies $\tilde s(\bar e) = \overline{s(e)}$ for all $e \in E(\Lambda_n)$.

Assume that $\tilde s(\bar e) = \tilde s(\bar f)$ for some $e, f \in E(\Lambda_n)$ where  $e = \{u,v\}$ and  $f = \{x,y\}$. Let $i,j \in [n]$ be such that $u_i \ne v_i$ and $x_j\ne y_j$. Then $\{u_i,v_i\} = \{x_j,y_j\} = \{0,1\}$ and $u_{i-1} = v_{i-1} = u_{i+1} = v_{i+1} = x_{j-1} = y_{j-1} = x_{j+1} = y_{j+1} = 0$. Without loss of generality assume $u_i = x_j = 1$.  Since  $\overline{s(e)} = \overline{s(f)}$, we distinguish two cases:

{\em Case 1:} If $s(e) = s(f)$ then $ u_{i+2}u_{i+3} \cdots u_{i+n-2} =  x_{j+2}x_{j+3} \cdots x_{j+n-2}$, hence
\begin{eqnarray*}
\alpha^{2-i}(u) &=&  u_{i-1}u_i u_{i+1} u_{i+2} \cdots u_{i+n-2}\ =\ 010 u_{i+2} \cdots u_{i+n-2} \\
&=&  x_{j-1} x_j x_{j+1} x_{j+2}x_{j+3} \cdots x_{j+n-2}\ =\ \alpha^{2-j}(x),
\end{eqnarray*}
and similarly, $\alpha^{2-i}(v) = \alpha^{2-j}(y)$. Consequently,  $\alpha^{2-i}(e) = \alpha^{2-j}(f)$.

{\em Case 2:} If $s(e) = \beta s(f)$ then $ u_{i+2}u_{i+3} \cdots u_{i+n-2} =  x_{j+n-2} \cdots x_{j+3} x_{j+2}$, hence
\begin{eqnarray*}
\alpha^{2-i}(u) &=&  u_{i-1}u_i u_{i+1} u_{i+2} \cdots u_{i+n-2}\ =\ 010 u_{i+2} \cdots u_{i+n-2} \\
&=&  x_{j+1} x_j x_{j-1} x_{j+n-2}\cdots x_{j+3}x_{j+2} \\ 
&=& \beta( x_{j+2}x_{j+3}\cdots  x_{j+n-2} x_{j-1} x_j x_{j+1})\ =\ \beta\alpha^{-j-1}(x)
\end{eqnarray*}
by (\ref{alphaj}), and similarly, $\alpha^{2-i}(v) = \beta\alpha^{-j-1}(y)$. Consequently,  $\alpha^{2-i}(e) = \beta\alpha^{-j-1}(f)$.

In either case, $\tilde s(\bar e) = \tilde s(\bar f)$ implies that $\bar e = \bar f$, hence $\tilde s$ is injective.

Now let $u \in V(\Gamma_{n-3})$. Then $e = \{010u,000u\} \in E(\Lambda_n)$ and $s(e) = u$, implying that $s$ and  $\tilde s$ are surjective. Thus $\tilde s$ is bijective and  $o_E(\Lambda_n) = o_V(\Gamma_{n-3})$.

Finally, (ii) follows from (i) and Theorem \ref{thm:number-of-orbits-Gamma} for $n \ge 3$, and can be verified directly for $n \le 2$.
\qed

\begin{theorem}
\label{LucasEorbitsize}
For all $n \in {\mathbb N}$,
\[
\{|X|;\ X \in {\cal O}_E(\Lambda_n)\}\ \subseteq\ \{n,2n\},
\]
with equality if and only if $n\ge 5$.
\end{theorem}

\proof
For $n \le 4$ this can be checked directly. Now assume that $n \ge 5$ and let $e=\{u,v\}\in E(\Lambda_n)$. Then by~\cite[Corollary~2~(ii)]{pasa-2002}, at least one of the strings $u$ and $v$ is primitive, say $u$. By Proposition~\ref{size2cu}, $|\bar u| = p(u)$ or $|\bar u| = 2p(u)$. Moreover, because $u$ is primitive,  $p(u) = |u| = n$. We conclude that $|\bar u| \in \{n, 2n\}$.

Denote by $S_u$ resp.\ $S_e$ the stabilizer of $u$ resp.\ $e$ under the action of $D_n$ on $V(\Lambda_n)$ resp.\ 
$E(\Lambda_n)$. Let $g \in S_e$. Then $\{g(u),g(v)\} = \{u,v\}$. Since $w(u)\ne w(v)$ and $g$ preserves weight, it follows that $g(u) = u$, so $g \in S_u$. Consequently $S_e \subseteq S_u$, hence $S_e \le S_u$ and $|S_e|$ divides $|S_u|$. From $|\bar e|\,|S_e| = |D_n| = |\bar u|\,|S_u|$ it now follows that $|\bar u|$ divides $|\bar e|$, so $n$ divides $|\bar e|$ as well. But $|\bar e|$ divides $2n$, hence $|\bar e| \in \{n, 2n\}$.  

To see that both $n$ and $2n$ indeed appear as orbit sizes, consider the edges $e=\{0^n, 10^{n-1}\}$ and $f=\{10^{n-1},1010^{n-3}\}$. Then it is straightforward to check that $|\bar e| = n$ and $|\bar f| = 2n$. 
\qed

\begin{corollary}
For all $n\in {\mathbb N}$,
\[
\begin{array}{lll}
o_E(\Lambda_n, n) &=&  F_{\lfloor \frac{n+1+(-1)^n}{2}\rfloor}, \\ 
o_E(\Lambda_n, 2n) &=& \frac{1}{2} \left(F_{n-1} - F_{\lfloor \frac{n+1+(-1)^n}{2}\rfloor}\right). \\ 
\end{array}
\]
\end{corollary}

\proof
This follows from (\ref{abE})  which in this case transforms into
\[
\begin{array}{lll}
n\, o_E(\Lambda_n, n) + 2n\, o_E(\Lambda_n, 2n) &=& n F_{n-1}, \\
o_E(\Lambda_n, n) + o_E(\Lambda_n, 2n) &=& \frac{1}{2} \left(F_{n-1} + F_{\lfloor \frac{n+1+(-1)^n}{2}\rfloor}\right). \qquad \qquad \qquad \mathqed
\end{array}
\]
\begin{table}[h]
\[
\begin{array}{c|rrrrrrrrrrrrrrrr}
n & 1& 2& 3& 4& 5& 6& 7& 8& 9& 10& 11& 12& 13& 14& 15& 16 \\
\hline
o_E(\Lambda_n) & 0& 1& 1& 2& 2& 4& 5& 9& 12& 21& 30& 51& 76& 127& 195& 322 \\
o_E(\Lambda_n, n) & 0& 1& 1& 2& 1& 3& 2& 5& 3& 8& 5& 13& 8& 21& 13& 34 \\
o_E(\Lambda_n, 2n) & 0& 0& 0& 0& 1& 1& 3& 4& 9& 13& 25& 38& 68& 106& 182& 288
\end{array}
\]
\caption{The numbers of all orbits, orbits of size $n$, and orbits of size $2n$ in $E(\Lambda_n)$ for $n \le16$ }
\end{table}

\section*{Acknowledgements} 

The research of Ali Reza Ashrafi and Khadijeh Fathalikhani was partially supported by the University of Kashan under the grant no 364988/9. The research of Sandi Klav\v zar and Marko Petkov\v sek was partially supported by the Ministry of Science of Slovenia under the grants P1-0297 and P1-0294.

\end{document}